\documentclass[a4paper, 12pt, fleqn]{article}

\usepackage[fleqn]{amsmath}
\usepackage{alltt}
\usepackage{amssymb}

\newtheorem{thm}{Theorem}
\newtheorem{lemma}[thm]{Lemma}

\begin{document}

\title{Weights in Codes and Genus 2 Curves}

\author{Gary McGuire\\
Department of Mathematics\\
NUI Maynooth\\
Co. Kildare\\
Ireland
\and
Jos\'e Felipe Voloch\\
Department of Mathematics\\
University of Texas\\
Austin, TX 78712\\
USA}

\maketitle
\begin{abstract}
\noindent
We discuss a class of binary cyclic codes and their dual codes.
The minimum distance is determined using algebraic geometry,
and an application of Weil's theorem.
We relate the weights
appearing in the dual codes to the number of rational points
on a family of genus 2 curves over a finite field.
\end{abstract}

\newpage

\section{Introduction}
Let $\mathbb{F}_q$ be a finite field with $q$ elements.
In this article, $q$ will be a power of $2$, say $q=2^m$,
and $\alpha$ will be a generator for the multiplicative
group $\mathbb{F}_q^*$.
Let $m_i(x)$ denote the minimal polynomial of $\alpha^i$
over $\mathbb{F}_2$.  
Cyclic codes of length $n$ are ideals in $\mathbb{F}_2 [x]/(x^n-1)$.
We use the natural basis $1, x, x^2, \ldots ,x^{n-1}$, and
we sometimes identify a polynomial 
$c_0+c_1x+\cdots +c_{n-1}x^{n-1}$
with the vector $(c_0, c_1, \ldots ,c_{n-1})$.
We label the coordinates by the elements of $\mathbb{F}_q^*$.

The cyclic code of length $2^m-1$ generated by $m_1(x)$
is called the (binary) Hamming code.
The cyclic code $B=B_m$ of length $2^m-1$ generated by $m_1(x)m_3(x)$
is called the 2-error-correcting BCH code.
The weights appearing in the dual code $B_m^\perp$ were
determined by Kasami \cite{K}.
There are exactly three nonzero weights when $m$ is odd,
and five weights when $m$ is even.
The cyclic code $M=M_m$ of length $2^m-1$ generated by $m_1(x)m_{-1}(x)$
is known as the Melas code.
The weights appearing in $M_m^\perp$ were determined by
Lachaud and Wolfmann \cite{LW} using results on elliptic curves.
In contrast to $B_m^\perp$, there are many weights in $M_m^\perp$.
Indeed, all even numbers between $q/2-\sqrt{q}+1/2$ and 
$q/2+\sqrt{q}+1/2$ occur.
A uniform treatment of these codes was given by Schoof \cite{S}.
In his paper Schoof says ``It would be very interesting to extend the
methods of this paper to other families of cyclic codes.
This seems difficult since it involves, in general,
curves of genus larger than 1 \ldots''

In this article we consider the 
cyclic code $C=C_m=B_m\cap M_m$, which has
length $2^m-1$ and is generated by 
$m_1(x)m_{-1}(x)m_3(x)$. We assume $m>2$ to ensure that the three
factors of the generator polynomial of $C$ are distinct.
We will discuss the minimum distance of $C$ in section \ref{mind},
using algebraic geometry.
In sections \ref{dual} and \ref{dual2}
we will determine the weights appearing
in the dual code $C^\perp$, by relating the weights to
curves of genus $2$, realising the suggestion of Schoof in the above quote.
For $m$ even we have a precise description of all the weights,
but not for $m$ odd.
The next step would be to compute the weight distributions of these
codes but this appears to be quite difficult.

\section{The Minimum Distance of the Codes $C$}\label{mind}

In this section we investigate the minimum distance of $C$.
We show below that $B$ has minimum distance 5, and it
is not hard to show that $M$ has minimum distance $5$
when $m$ is odd, and distance 3 when $m$ is even.
Since $C=B\cap M$,
one might hope that $C$ has distance $7$, at least when $m$ is odd.
However, we will show that the minimum distance of $C$ is $5$
for $m\geq 16$. 

The computer algebra package \emph{Magma} shows that
$C$ has minimum distance $7$ when $m=6,7$, but that
$d(C)=5$ when $m=5,8,9$.
We presume that $d(C)=5$ when $10\leq m \leq 15$
but we have not checked this.

The roots of the generator polynomial of a cyclic code
are called the \emph{zeros} of the code.
Determining the minimum distance of a cyclic code from
its zeros is very difficult in general. 
One result on this problem is known as the BCH bound,
see \cite{MS} for example.
We use $wt(c)$ to denote the weight of a codeword $c(x)$.

\begin{thm} (BCH bound)
Let $f(x)$ be a codeword in a binary cyclic code of length $n=2^m-1$.
If $s$ consecutive powers of $\alpha$ are roots of $f$,
then $wt(f)\geq s+1$.
\end{thm}

It follows from the BCH bound that the 2-error-correcting
BCH code $B_m$ has $d\geq 5$,
since $\alpha, \alpha^2, \alpha^3, \alpha^4$ are roots of
$m_1(x) m_3(x)$.
Since $C \subseteq B$, $d(C) \geq 5$.

A codeword of even weight in $C$ has amongst its roots
$\alpha^j$ for $j=-2, -1$, $0, 1, 2$, $3, 4$.
By the BCH bound this codeword has weight $\geq 8$.
Thus there are no codewords of weight $6$ in $C$.
We now study codewords of weight $5$.

We define the polynomials
\[
f(x,y,z)=x+y+z+x^2+y^2+z^2+x^2y+x^2z+y^2x+y^2z+z^2x+z^2y
\]
and
\[
g(x,y,z)=x^2y+x^2z+y^2x+y^2z+z^2x+z^2y+xyz+xy+xz
\]
\[
+yz+x^2yz+xy^2z+xyz^2
\]
over a field of characteristic $2$.
Let $K$ be the algebraic closure of $\mathbb{F}_2$.
We define an algebraic curve $X$ by
\[
X=\{(x,y,z)\in K^3 : f(x,y,z)=0 \textrm{  and  } g(x,y,z)=0\}
\]
and define $X_m$ to be the set of points on $X$ that have
coordinates in $\mathbb{F}_{2^m}$.

\begin{lemma}\label{pclemma}
The cyclic code $C$ of length $2^m-1$ has minimum distance $5$
if and only if there are rational points $(x,y,z)$ on $X_m$
with the property that $0, 1, x, y, z, 1+x+y+z$ are pairwise distinct.
\end{lemma}

Proof:
A parity check matrix for $C$ is
\[  \left[%
\begin{array}{ccccccc}
  1 & \alpha & \alpha^2 & \cdots &  \alpha^i & \cdots & \alpha^{2^m-2} \\
  1 & \alpha^{3} & \alpha^{6} & \cdots & \alpha^{3i} & \cdots & \alpha^{3(2^m-2)} \\
  1 & \alpha^{-1} & \alpha^{-2} & \cdots & \alpha^{-i} & \cdots & \alpha^{-(2^m-2)} \\
\end{array}%
\right]
\]
and it follows that codewords of weight $5$ 
with a 1 in position 1
are in one-to-one correspondence 
with field elements $x, y, z, w \in \mathbb{F}_{2^m}$ such that
\begin{eqnarray}
1+x+y+z+w&=&0\label{first}\\
1+x^3+y^3+z^3+w^3&=&0\label{second}\\
1+x^{-1}+y^{-1}+z^{-1}+w^{-1}&=&0\label{third}
\end{eqnarray}
and $0, 1, x, y, z, w$ are pairwise distinct.

Substituting $1+x+y+z$ for $w$ in equation (\ref{second}) gives
\[
1+x^3+y^3+z^3+(1+x+y+z)^3=0
\]
or
\[
x+y+z+x^2+y^2+z^2+x^2y+x^2z+y^2x+y^2z+z^2x+z^2y=0
\]
which leads to the definition of $f$.

Multiplying (\ref{third}) by $xyzw$ gives
\[
xyzw+yzw+xzw+xyw+xyz=0
\]
and substituting for $w$ now gives
\[
(1+x+y+z)(xyz+yz+xz+xy)+xyz=0.
\]
Expanding this out leads to
\[
x^2y+x^2z+y^2x+y^2z+z^2x+z^2y+xyz+xy+xz+yz+x^2yz+xy^2z+xyz^2=0
\]
which is where we obtain the definition of $g$.

The proof is complete when we observe that the steps in deriving
$f$ and $g$ are reversible; given a point on $X_m$ with
$0, 1, x, y, z, 1+x+y+z$  distinct, one can recover a codeword of
weight 5 with a 1 in position 1.
Since $C$ is cyclic, any weight 5 codeword has a cyclic shift
with a 1 in position 1.
\hfill $\Box$

\bigskip

We will apply Weil's theorem to $X$.
Normally Weil's theorem is applied to nonsingular curves,
but a straightforward check via the Jacobian matrix shows
that $X$ has exactly four singular points.
However, the nonsingularity hypothesis
in Weil's theorem can be replaced by
absolute irreducibility, and we show next that this
indeed holds for our curve $X$.

\begin{lemma}\label{absred}
The curve $X$ is absolutely irreducible.
\end{lemma}

Proof:
Define 
\[
a(x,y)=1+x+y, \qquad c(x,y)=xy+x+y,
\]
and
\[
h(x,y)=(y^2+y+1)x^3+(y^3+1)x^2+(y^3+y)x+(y^3+y^2).
\]
With $f$ and $g$ as above, we verify that
$ag+cf=h$,
which is independent of $z$.
It is straightforward to check that $h$ is absolutely irreducible.
(This can be done by hand or using a computer package such as \emph{Magma}.
Since $h$ is of degree 3 in $x$ it is enough to check irreducibility over
$\mathbb{F}_8$.   \emph{Magma} also shows that $h=0$ has genus 3.)


Let $Y$ be the plane curve $h=0$.
Since $f=az^2+a^2z+b$ and $g=cz^2+acz+d$ for some polynomials
$b(x,y)$ and $d(x,y)$,
projection on the $x,y$ plane gives a map from $X$ to $Y$, which is of
degree 2. Since we already know that $Y$ is absolutely irreducible,
we get that either $X$ is also absolutely irreducible or it has two
components.

Let $w$ be a primitive 3rd root of $1$ in $GF(4)$. Then $h(w,w^2)=0$
while $h_x(w,w^2)=w$ and $h_y(w,w^2)=w^2$. So the point $(w,w^2)$ is a
smooth point of $Y$ with tangent $y=w^2x+w$.

In the equation $f=0$ make the substitution $v=z/a$, and the equation 
becomes $v^2+v=b/a^3$.

Consider $b/a^3$ as a function on $Y$, and consider its
behaviour near the point $(w,w^2)$.  Note that $a$ vanishes at $(w,w^2)$
but since $a=0$ is not the tangent to $Y$ at $(w,w^2)$, the function $a$
has a simple zero at $(w,w^2)$. On the other hand $b(w,w^2)=1$,
so $b/a^3$ has a triple pole at $(w,w^2)$. 
However, if $v^2+v$ has a pole at a point $P$ then the pole
must have even order (the order is $2t$, where $t$ is the order
of $v$ at $P$).
Thus there cannot
be a function $v$ on $Y$ with $v^2+v=b/a^3$.
This means that the polynomial
$v^2+v+b/a^3$ is irreducible over the function field of $Y$, which entails
that $X$ is absolutely irreducible.
\hfill $\Box$

\bigskip

\begin{thm}
The cyclic code $C$ of length $2^m-1$ has minimum distance $5$
for all $m\geq 16$.
\end{thm}

Proof:
By Lemma \ref{pclemma} we must show that $X_m$ has points $(x,y,z)$
with $0, 1, x, y, z, 1+x+y+z$  distinct, for all $m$ sufficiently
large.  By Lemma \ref{absred} we know that $X_m$ is absolutely
irreducible.  
Let $N_m=|X_m|$.
By Weil's theorem,
\[ 
|N_m - (2^m+1)|  \leq 2g \sqrt{2^m} + C
\]
where $g$ is the genus of (a smooth model of) $X$ and $C$ is a constant
independent of $m$ which can be given in terms of the degree of $X$.

The number of points on $X_m$ such that
$0, 1, x, y, z, 1+x+y+z$  are not distinct is 4.
This is straightforward to check using such 
factorizations as
$f(0,y,z)=(y+1)(z+1)(y+z)$, and we omit the details
(the four points are $(0,0,0)$, $(1,0,0)$,
$(0,1,0)$, and $(0,0,1)$).

It follows from the previous two paragraphs that there
are points on $X_m$ with 
$0, 1, x, y, z, 1+x+y+z$ distinct for all
$m$ sufficiently large.

Using a refined version of Weil's theorem from \cite{AP},
we obtain $|N_m - (2^m+1)|  \leq 220\sqrt{2^m}$.
It follows from
this inequality that $N_m>4$ once $m\geq 16$.
\hfill $\Box$

\bigskip

It can be easily shown that the genus of $X$ is between $11$ and $13$,
but we have not computed its exact value.

\vskip1in

\section{The Weights in the Dual Codes $C^\perp$, $m$ even}\label{dual}

Let $q=2^m$.  By Delsarte's theorem (see \cite{S} or \cite{MS}), 
\[
C^\perp = \{(Tr(a/x+bx+cx^3))_{x\in \mathbb{F}_q^*} : 
a,b,c\in \mathbb{F}_{q}\}.
\]
Knowing the weights in $C^\perp$ is equivalent to knowing
how many $0$'s are in a typical codeword.
By Hilbert's Theorem 90, we want to know how many solutions
there are to
\begin{equation}
y^2+y=\frac{a}{x} +bx+cx^3\label{fourth}
\end{equation}
over $\mathbb{F}_{2^m}$. If we denote by $N$ the number of rational points
in a complete smooth model of the above curve then
the weight of the vector whose entries are $Tr(a/x+bx+cx^3)$ as
we vary $x\in \mathbb{F}_q^*$, is $q-1-(N-2)/2=q-N/2$.

Recall that every curve has an abelian variety associated
to it called its Jacobian.
An abelian variety $A$ over a field of characteristic $p>0$
is said to have $p$-rank $s$
if the subgroup of points of order $p$ of $A$ (over an algebraically
closed field of definition) has cardinality $p^s$.
By the two-rank of a curve we mean the two-rank of its
Jacobian. 

\begin{lemma}\label{char}
Curves of the form $(\ref{fourth})$ can be characterised as curves defined
over $\mathbb{F}_{2^m}$ of genus $2$, two-rank $1$, whose 
number of rational points is divisible by $4$. 
\end{lemma}

Proof: From \cite{CNP}, it follows that 
a curve of genus 2 and two-rank 1
has an equation $y^2+y=a/x+bx+cx^3+d$.
Let us now show that we may take $d=0$ when $N \equiv 0 \mod 4$.
As the number of points
is zero modulo $4$, $\sum_{x\in \mathbb{F}_q^*} Tr(a/x+bx+cx^3+d) = 0$. But
$\sum_{x\in \mathbb{F}_q^*} 1/x = \sum_{x\in \mathbb{F}_q^*} x 
= \sum_{x\in \mathbb{F}_q^*} x^3=0$,
if $q>4$ so we get $Tr(d)=0$ and $d=e^2+e$ for some $e$ and a change
of variable $y \mapsto y+e$ puts the equation in the form stated with $d=0$.
Conversely, if $d=0$ then $N \equiv 0 \pmod4$ as
$\sum_{x\in \mathbb{F}_q^*} Tr(a/x+bx+cx^3) = 0$.
\hfill $\Box$

\bigskip

An abelian variety is called simple if is not isogenous to a product of
abelian varieties of smaller dimension. Maisner and Nart classified
which isogeny classes of simple abelian surfaces of $p$-rank one
contain Jacobians.

\begin{thm}\label{mn}
(Maisner-Nart, \cite{MN}) Let $q=2^m$.
There exists a curve of the form $(\ref{fourth})$ with $N=q+1+a_1$ 
points over $\mathbb{F}_{2^m}$ and simple Jacobian if and only if 
\begin{enumerate}
\item $a_1$ is odd
\item $|a_1| \le 4\sqrt{q}$ 
\item there exists an integer $a_2$
such that 
\begin{enumerate}
\item $2|a_1|\sqrt{q}-2q \le a_2 \le a_1^2/4+2q$ 
\item $a_2$ is divisible by
$2^{\lceil m/2 \rceil}$ 
\item $\Delta=a_1^2-4a_2+8q$ is not a square in $\mathbb{Z}$
\item $\delta=(a_2+2q)^2-4qa_1^2$ is not a square in $\mathbb{Z}_2$.
\end{enumerate}
\end{enumerate}
\end{thm}

This statement combines Lemma 2.1, Theorem 2.9 part (M) 
and Corollary 2.17 of  \cite{MN}, and our Lemma \ref{char}.

\begin{lemma}\label{J}
Let $q=2^m$ where $m$ is even.
Then each even number in the interval
$[q/2-2\sqrt{q}+q^{\frac{1}{4}}-1/2,q/2+2\sqrt{q}-q^{\frac{1}{4}}-1/2]$
occurs as a weight in $C^\perp$, and these weights arise
from curves of type $(\ref{fourth})$ whose Jacobian is simple.
\end{lemma}

Proof:
Assume  that $m$ is even.
If $a_1,a_2$ satisfy the conditions of theorem \ref{mn} then
$|a_1| \le 4\sqrt{q}-2q^{1/4}$. Indeed, if $a_2=2|a_1|\sqrt{q}-2q$ then
$\Delta=(|a_1|+4\sqrt{q})^2$ is a square, so it is ruled out. Thus
$2|a_1|\sqrt{q}-2q+\sqrt{q} \le a_2 \le a_1^2/4+2q$, which leads to the
stated inequality.

Conversely, if $a_1$ satisfies
the inequality $|a_1| \le 4\sqrt{q}-2q^{1/4}$,
let $a_2=2|a_1|\sqrt{q}-2q+\sqrt{q}$.
We must check that $\Delta$ and $\delta$ are not squares
in $\mathbb{Z}$ and $\mathbb{Z}_2$ respectively.

First, substitution gives $\Delta=(4\sqrt{q}-|a_1|)^2-4\sqrt{q}$.
Suppose $\Delta=t^2$ where $t$ is a positive integer.
Then $(4\sqrt{q}-|a_1|-t)(4\sqrt{q}-|a_1|+t)=4\sqrt{q}$.
By unique factorization in $\mathbb{Z}$, we conclude
$4\sqrt{q}-|a_1|-t=2^k$ and $4\sqrt{q}-|a_1|+t=2^\ell$
for some positive integers $k$ and $\ell$ with $k+\ell=2+m/2$.
Adding gives $2(4\sqrt{q}-|a_1|)=2^k+2^\ell$.
Since $a_1$ is odd, one of $k$ and $\ell$ must be 1.
If $\ell=1$, then $4\sqrt{q}-|a_1|+t=2$,
so $4\sqrt{q}-|a_1|=t=1$, a contradiction.
Suppose now that $k=1$ (so $\ell=1+m/2$).
It follows that $t=4\sqrt{q}-|a_1|-2$.
Substituting this value for $t$ into $4\sqrt{q}-|a_1|+t=2^\ell$
yields $|a_1|=3\sqrt{q}-1$.
Thus, if $|a_1|\not=3\sqrt{q}-1$ we have shown that
$\Delta$ is not a square.

If $|a_1|=3\sqrt{q}-1$ then choose
$a_2=2|a_1|\sqrt{q}-2q+2\sqrt{q}$.
A similar argument as above leads to a contradiction.

Next, substituting for $a_2$ gives
\[
\delta =(a_2+2q)^2-4qa_1=(2|a_1|\sqrt{q}+\sqrt{q})^2-4qa_1^2=
q(1+4|a_1|).
\]
It is well known (see \cite{serre} ch.\ II for example)
that an element $2^ru$ (where $u$ is a unit)
of $\mathbb{Z}_2$ is a square if and
only if $r$ is even and $u\equiv 1\pmod8$.
Since $a_1$ is odd it follows that $\delta$ is not a square.
A similar argument works in the case 
$a_2=2|a_1|\sqrt{q}-2q+2\sqrt{q}$.

\hfill $\Box$

\bigskip

We still need to analyse which weights come from curves whose Jacobian
is non-simple. 
We do this in the proof of theorem \ref{meven}.
We note that by \cite{MN} corollary 2.17 the field
of definition does not matter to determine simplicity.

\begin{thm}\label{meven}
Let $q=2^m$ where $m$ is even,
let $I=[q/2-2\sqrt{q},q/2+2\sqrt{q}-1]$ and
$J=[q/2-2\sqrt{q}+q^{\frac{1}{4}}-\frac{1}{2},
q/2+2\sqrt{q}-q^{\frac{1}{4}}-\frac{1}{2}]$.
Then all weights in $C^\perp$ are even integers in $I$.
All even integers in $J$ do occur as weights, and
an even integer in $I\setminus J$ occurs as a weight if and only if
it has the form $q/2+(\pm 2\sqrt{q}+a+1)/2$ where $a\equiv 3\pmod4$
and $\pm 2\sqrt{q} -a$ is not squarefree.
\end{thm}

Proof:
We continue the notation from before.
The weights are the numbers $q-N/2$, where
$N=q+1+a_1$ ranges over the number
of points on curves of type (\ref{fourth}). From theorem \ref{mn} $a_1$ 
is odd and $|a_1|\leq 4\sqrt{q}$.
Thus $-4\sqrt{q}+1\leq a_1\leq 4\sqrt{q}-1$, and this is equivalent
to saying that the weights lie in $I$.
All weights in $C^\perp$ are even since $1$ is a zero of the code.
(This entails $N\equiv 0\pmod4$, which means $a_1\equiv 3 \pmod4$.)

By Lemma \ref{J} all weights in $J$ do occur as weights.

We now study curves of type (\ref{fourth}) whose Jacobian
is not simple.  In this case the Jacobian must be isogenous to
$E'\times E$, where $E'$ is an elliptic curve of two-rank 0
(a supersingular elliptic curve) and $E$ is an elliptic curve of
two-rank 1 (an ordinary elliptic curve).
It is known (see \cite{S} for example) that a 
supersingular elliptic curve $E'$
has $q+1-a'$ points, where $a'\in \{0,\pm \sqrt{q}, \pm 2\sqrt{q}\}$
(as $m$ is even).  It is also known by results of Honda and Tate
that an ordinary elliptic curve $E$ exists with $q+1-a$ points whenever
$a$ is odd and $|a|\leq 2\sqrt{q}$.
We will investigate when we can
construct a curve of genus $2$ having $N=q+1-a'-a$ points over 
$\mathbb{F}_q$ 
whose Jacobian is isogenous to $E'\times E$. To do this
we apply the construction of \cite{FK}, section 1. There it is proved
that such a curve exists if and only if,
for some odd prime $p$, there is an isomorphism of Galois modules
between $E'[p]$ and $E[p]$ reversing the Weil pairing. 

We will restrict ourselves to the case that $E'$ has $q+1\pm 2\sqrt{q}$ points.
In the other cases the construction can be done whenever $a-a' \ne \pm 1$
but it leads to weights in the interval $J$, which are therefore not 
interesting for our purposes. Returning to the case in question,
the action of Frobenius on $E'[p]$ is multiplication
by a scalar $k=\pm\sqrt{q}$. We need to have the same be true of
$E[p]$, and then any isomorphism of groups
between $E'[p]$ and $E[p]$ reversing the Weil pairing will automatically
preserve Galois action and we will be done. 

To have Frobenius on $E[p]$
be multiplication by $k$ we must have $k^2-ak+q \equiv 0 \mod p^2$.
Indeed, $k^2-ak+q$ is the number of points in the kernel of $\pi - k$
on $E$, where $\pi$ denotes the Frobenius automorphism on $E$.
Conversely, if $k^2-ak+q \equiv 0 \mod p^2$, we will show that either
Frobenius on $E[p]$ is multiplication by $k$ or there is an elliptic
curve isogenous to $E$ with this property. The congruence implies that
the characteristic polynomial of $\pi$ on $E[p]$ is $(x-k)^2$.
Assume now that Frobenius on $E[p]$ is not multiplication by $k$. Then
$\pi - k$ has a kernel $\Gamma$ on $E[p]$ which is also the image of
$\pi - k$ on $E[p]$. Thus $\Gamma$ is invariant under $\pi$ and hence
$\overline{E}=E/\Gamma$ is defined over 
$\mathbb{F}_q$ and is isogenous to $E$. Now
$\pi - k=0$ on $E[p]/\Gamma \subset \overline{E}[p]$ 
and by the same argument as above
(since the congruence holds modulo $p^2$)
$\pi - k=0$ on a cyclic subgroup of $E[p^2]$ which projects to a
different subgroup of $\overline{E}[p]$, 
thus $\pi - k=0$ on $\overline{E}[p]$.
To summarize, we can construct the curve of genus two if
$a'\equiv a \mod p^2$ for some prime $p$, when $a' = \pm 2\sqrt{q}$.

Therefore a value of $a_1=\pm2\sqrt{q}+a$ is realisable
from this construction if and only if $\pm2\sqrt{q}-a$ is not squarefree.
\hfill $\Box$

\bigskip

Here are the lists of weights in a few cases.
\bigskip

\begin{tabular}{|c|c|c|c|c|} 
\hline
$q$ & $2^6$ & $2^8$ & $2^{10}$ & $2^{12}$ \\
\hline
$I$& $[16,47]$ & $[96,159]$ & $[448,575]$&$[1920,2175]$ \\
\hline
$J$ & $[19,44]$ & $[100,155]$ & $[454,569]$&$[1928,2167]$ \\
\hline
weights in $I\setminus J$ & none &none &452 &1924 \\
\hline
\end{tabular}

\bigskip

We point out that the weights are not necessarily all
the even numbers in an interval, as illustrated by
the $q=2^{12}$ case.

\section{The Weights in the Dual Codes $C^\perp$, $m$ odd}\label{dual2}

Let us consider now the case $m$ odd.

\begin{thm}\label{modd}
Let $q=2^m$ where $m$ is odd,
let $I=[q/2-\lfloor 2\sqrt{q}\rfloor,q/2+\lfloor 2\sqrt{q}\rfloor-1]$ and
$J=[q/2-2\sqrt{q}+(8q)^{\frac{1}{4}}-\frac{1}{2},
q/2+ 2\sqrt{q}-(8q)^{\frac{1}{4}}-{\frac{1}{2}}]$.
Then all weights in $C^\perp$ are even integers in $I$, and
all even integers in $J$ do occur as weights.
\end{thm}

Proof:
We need
only to consider the values of $a_1$ afforded by Theorem \ref{mn}, since
the curves with split Jacobian will have number of points of the form
$q+1+a,q+1\pm \sqrt{2q} +a$, for some $a$ satisfying $|a| \le 2\sqrt{q}$
which will lead to weights in $J$. Note also that
we can improve the inequality on $a_1$ to 
$|a_1|\leq 2\lfloor 2\sqrt{q}\rfloor$,
as noted in \cite{MN}. This leads to interval $I$ and the first statement
of the Theorem.

Let $q' = 2^{(m+1)/2}$. Let $a_1$ be an odd integer and let
$a_2$ be any integer divisible by $q'$, and put 
$\delta = (a_2 + 2q)^2-4qa_1^2$. We will show that 
$\delta$ is not a square in $\mathbb{Z}_2$. Consider first the case
where $a_2=q'u$, $u$ odd. Recall that $(q')^2 = 2q$.
Then $\delta/2q \equiv u^2-2a_1^2 \mod 8$.  As $u^2=a_1^2 \equiv 1 \mod 8$
we get $\delta/2q \equiv 7 \mod 8$ and $\delta$ 
is not a square in $\mathbb{Z}_2$.
If $a_2=2q'u$, $u$ odd, then $\delta/4q = 2u^2-a_1^2$ is odd.
Again $\delta$ is not a square in $\mathbb{Z}_2$ as
$\delta = 2^ru$ where $r$ is odd and $u$ is a unit. Finally, 
if $a_2/q' \equiv 0 \mod 4$ then $\delta/2q \equiv -2a_1^2 \equiv -2 \mod 8$, 
so $\delta$ is again not a square in $\mathbb{Z}_2$.

If, now we assume further that $a_1 \in J$ then there exists an 
integer $a_2$ such that 
$a_2,a_2+q'$ satisfy conditions (a) and (b) of Theorem \ref{mn}.
By the above argument they also satisfy condition (d). 
We will show that at least one of them satisfies condition (c).

Suppose neither of them satisfies condition (c).
Let $\Delta(b)=a_1^2-4b+8q$.
If $\Delta(a_2)=u^2,\Delta(a_2+q')=v^2$ for positive integers $u,v$
then $u^2-v^2=4q'$. It follows that $u-v=2^r,u+v=2^s$ 
for some integers $r,s$, where $r+s=(m+5)/2$.
So $v=2^{s-1}-2^{r-1}$. However, since $a_1$
is odd, it follows that $v^2=\Delta(a_2+q')$ is also odd, so $v$ is odd
and thus $r=1$ and $s=(m+3)/2$ which implies that $u=q'+1$ and so
$a_1^2 \equiv \Delta(a_2)=u^2 \equiv 1 \mod 2q'$.
Since $a_1 \equiv 3 \mod 4$ it follows that $a_1 \equiv -1 \mod q'$.
On the other hand $|a_1| \le 4\sqrt{q}= 2\sqrt{2}q'$, and we conclude that
$a_1 = -1+nq'$, where $n\in\{ 0,\pm 1,\pm 2\}$. 
We can then conclude that 
there exists a possibly different integer $a_2$ such that
$a_2,a_2+q',a_2+2q'$ satisfy conditions (a) and (b) of Theorem \ref{mn}.
By the above argument they also satisfy condition (d).
We proceed to show that at least one of them also satisfy condition (c).
If none of them satisfies condition (c) then we can apply the above
argument to both pairs $a_2,a_2+q'$ and $a_2+q',a_2+2q'$, but $u,v$ were
uniquely determined in terms of $q'$ above so we cannot have two such pairs.
The Theorem now follows from Theorem \ref{mn}.
\hfill $\Box$

\bigskip

Here are the lists of weights in a few cases.
Again we note that the weights are not necessarily all
the even numbers in an interval, as illustrated by
the $q=2^{11}$ case.
\bigskip

\begin{tabular}{|c|c|c|c|c|} 
\hline
$q$ & $2^7$ &  $2^{9}$ & $2^{11}$ \\
\hline
$I$& $[42,85]$ &  $[211,300]$&$[934,1113]$ \\
\hline
$J$ & $[47,80]$ & $[219,292]$ & $[945,1102]$ \\
\hline
weights in $I\setminus J$ & 46, & 216,218,
 & 938,942,944,\\
&82,84& 294,296&1104,1106 \\
\hline
\end{tabular}

\bigskip

We do not have a precise description of the weights in $I\setminus J$,
unlike the $m$ even case.  
The entries in the table for $I\setminus J$ were 
determined by computer.


\bigskip

\begin{thebibliography}{99}



\bibitem{AP}
Y.\ Aubry and  M.\ Perret,
A Weil theorem for singular curves,
``Arithmetic, geometry and coding theory (Luminy, 1993),''  1--7,
de Gruyter, Berlin, 1996.

\bibitem{CNP}
G. Cardona, E. Nart and J. Pujolas
Curves of genus two over fields of even characteristic. preprint.
arXiv:math.NT/0210105 


\bibitem{FK}
G. Frey  and E. Kani, 
Curves of genus $2$ covering elliptic curves and an arithmetical
application.
in ``Arithmetic algebraic geometry (Texel, 1989),'' 153--176,
Progr. Math., 89, Birkh\"auser Boston, Boston, MA, 1991.

\bibitem{K}
T. Kasami,
Weight distributions of Bose--Chaudhuri--Hocquenghem codes,
in ``Combinatorial Mathematics and its Applications 
(Proc. Conf., Univ. North Carolina, Chapel Hill, N.C., 1967),'' 
335--357, Univ. North Carolina Press, Chapel Hill, N.C. 1969.

\bibitem{LW}
G.\ Lachaud and J. Wolfmann, The weights of the orthogonals of
the extended quadratic binary Goppa codes,
{\it IEEE Trans. Info. Th.} {\bf 36} No. 3 (1990) 686--692.

\bibitem{MN}
D. Maisner and E. Nart,
Abelian surfaces over finite fields as Jacobians,
{\it Experimental Math.} {\bf 11} No. 3 (2002) 321--338.

\bibitem{MS}
F.\ J.\ MacWilliams and N.\ J.\ A.\ Sloane, ``The Theory of
Error--Correcting Codes,'' North Holland, Amsterdam, 1977.

\bibitem{S}
R.\ Schoof, Families of curves and weight distributions of codes,
{\it Bull. AMS} {\bf 32} No. 2 (1995) 171--183.

\bibitem{serre} 
J.-P.\ Serre, ``A Course in Arithmetic,'' GTM 7, Springer-Verlag,
New York, 1973.

\end{thebibliography}
\end{document}